\newif\ifspringer
\springerfalse
\ifspringer

    \documentclass{svmult}

    \smartqed
    \usepackage{mathptmx}       % selects Times Roman as basic font
    \usepackage{helvet}         % selects Helvetica as sans-serif font
    \usepackage{courier}        % selects Courier as typewriter font
    \usepackage{type1cm}        % activate if the above 3 fonts are
    \renewcommand{\email}[1]{\emailname: #1} % change the email address font style
    
    \spdefaulttheorem{assumption}{Assumption}{\upshape \bfseries}{\itshape}
    \spdefaulttheorem{algo}{Algorithm}{\upshape \bfseries}{\itshape}

\else

    \documentclass{article}

    \usepackage[a4paper]{geometry}

    \newenvironment{acknowledgement}{\paragraph{Acknowledgement.}}{\par}
    \newcommand{\email}[1]{\texttt{#1}}

    \usepackage{amsthm}
    \theoremstyle{plain}

    \theoremstyle{definition}

    \theoremstyle{remark}

\fi

\usepackage{array,colortbl}
\usepackage{amsmath,amsfonts,amssymb,bm} % no amsthm, Springer defines Theorem, Lemma, etc themselves
\DeclareFontFamily{U}{mathx}{\hyphenchar\font45}
\DeclareFontShape{U}{mathx}{m}{n}{
      <5> <6> <7> <8> <9> <10>
      <10.95> <12> <14.4> <17.28> <20.74> <24.88>
      mathx10
      }{}
\DeclareSymbolFont{mathx}{U}{mathx}{m}{n}
\DeclareFontSubstitution{U}{mathx}{m}{n}
\DeclareMathAccent{\widecheck}      {0}{mathx}{"71}

\usepackage{makeidx}         % allows index generation
\usepackage{graphicx}        % standard LaTeX graphics tool
\usepackage[bottom]{footmisc}% places footnotes at page bottom

\usepackage{latexsym}
\usepackage{amsmath}
\usepackage{amsfonts}
\usepackage{amssymb}
\usepackage{bm}

\usepackage{url}
\usepackage{algorithm}
\usepackage{algorithmic}
\usepackage[misc,geometry]{ifsym}

\newcommand{\bsb}{{\boldsymbol{b}}}

\newcommand{\bse}{{\boldsymbol{e}}}

\newcommand{\bsm}{{\boldsymbol{m}}}

\newcommand{\bst}{{\boldsymbol{t}}}

\newcommand{\bsw}{{\boldsymbol{w}}}
\newcommand{\bsx}{{\boldsymbol{x}}}
\newcommand{\bsy}{{\boldsymbol{y}}}
\newcommand{\bsz}{{\boldsymbol{z}}}

\newcommand{\bszero}{{\boldsymbol{0}}} % vector of zeros
\newcommand{\bsbeta}{{\boldsymbol{\beta}}}
\newcommand{\bsgamma}{{\boldsymbol{\gamma}}}

\newcommand{\bsnu}{{\boldsymbol{\nu}}}

\newcommand{\bsDelta}{{\boldsymbol{\Delta}}}

\newcommand{\rd}{{\mathrm{d}}}

\newcommand{\bbN}{{\mathbb{N}}}

\newcommand{\bbP}{{\mathbb{P}}}

\newcommand{\bbR}{{\mathbb{R}}}

\newcommand{\bbZ}{{\mathbb{Z}}}
\DeclareSymbolFont{bbold}{U}{bbold}{m}{n}
\DeclareSymbolFontAlphabet{\mathbbold}{bbold}

\newcommand{\calB}{{\mathcal{B}}}

\newcommand{\calO}{{\mathcal{O}}}

\newcommand{\setu}{{\mathfrak{u}}}

\DeclareSymbolFont{bbold}{U}{bbold}{m}{n}
\DeclareSymbolFontAlphabet{\mathbbold}{bbold}

\usepackage{microtype} % good font tricks

\usepackage[colorlinks=true,linkcolor=black,citecolor=black,urlcolor=black]{hyperref}
\usepackage{bookmark}
\makeatletter % to avoid hyperref warnings:
  \providecommand*{\toclevel@author}{999}
  \providecommand*{\toclevel@title}{0}
\makeatother

\usepackage{enumerate} % By including this package, the enumeration labels do not go into the margin
\begin{document}

\ifspringer

    \title*{Application of Quasi-Monte Carlo Methods to PDEs with Random
    Coefficients -- an Overview and Tutorial}
    \titlerunning{QMC for PDEs with random coefficients -- an overview and
    tutorial}
    % Use \titlerunning{Short Title} for an abbreviated version of
    % your contribution title if the original one is too long
    \author{Frances Y.\ Kuo \and Dirk Nuyens}
    % Use \authorrunning{Short Title} for an abbreviated version of
    % your contribution title if the original one is too long
    \institute{
        Frances Y. Kuo (\Letter)
        \at School of Mathematics and Statistics,
        University of New South Wales, Sydney NSW 2052, Australia \\
        \email{f.kuo@unsw.edu.au}
        \and
        Dirk Nuyens
        \at Department of Computer Science, KU Leuven,
        Celestijnenlaan 200A, 3001 Leuven, Belgium \\
        \email{dirk.nuyens@cs.kuleuven.be}
    }

    \index{Kuo, Frances Y.}
    \index{Nuyens, Dirk}

\else

    \title{Application of Quasi-Monte Carlo Methods to PDEs with Random
    Coefficients -- an Overview and Tutorial}

    \author{Frances Y.\ Kuo\footnote{
            Frances Y. Kuo (\Letter):
            School of Mathematics and Statistics,
            University of New South Wales, Sydney NSW 2052, Australia,
            \email{f.kuo@unsw.edu.au}
        }
       \and Dirk Nuyens\footnote{
            Dirk Nuyens:
            Department of Computer Science, KU Leuven,
            Celestijnenlaan 200A, 3001 Leuven, Belgium,
            \email{dirk.nuyens@cs.kuleuven.be}
        }
    }

    \date{}

\fi

\maketitle

\abstract{This article provides a high-level overview of some recent works
on the application of quasi-Monte Carlo (QMC) methods to PDEs with random
coefficients. It is based on an in-depth survey of a similar title by the
same authors, with an accompanying software package which is also briefly
discussed here. Embedded in this article is a step-by-step tutorial of the
required analysis for the setting known as the uniform case with first
order QMC rules. The aim of this article is to provide an easy entry point
for QMC experts wanting to start research in this direction and for PDE
analysts and practitioners wanting to tap into contemporary QMC theory and
methods.}

\section{Introduction}

\emph{Uncertainty quantification} is the science of quantitative
characterization and reduction of uncertainties in both computational and
real world applications, and it is the source of many challenging high
dimensional integration and approximation problems. Often the high
dimensionality comes from uncertainty or randomness in the data, e.g., in
groundwater flow from permeability that is rapidly varying and uncertain,
or in financial mathematics from the rapid and often unpredictable changes
within markets. The input data may be a random variable or a random field,
in which case the derived quantity of interest will in general also be a
random variable or a random field. The computational goal is usually to
find the expected value or other statistics of these derived quantities.

A popular example is the flow of water through a disordered porous medium,
modeled by Darcy's law coupled with the mass conservation law, i.e.,
\begin{align*}
 q(\bsx,\omega) + a(\bsx,\omega)\,\nabla p(\bsx,\omega) &\,=\, 0\,, \\
 \nabla\cdot q(\bsx,\omega) &\,=\, 0\,,
\end{align*}
for $\bsx$ in a bounded domain $D\subset\bbR^d$, $d\in\{1,2,3\}$, and for
almost all events $\omega$ in the probability space
$(\Omega,\mathcal{A},\bbP)$. Here $q(\bsx,\omega)$ is the velocity (also
called the specific discharge) and $p(\bsx,\omega)$ is the residual
pressure, while $a(\bsx,\omega)$ is the permeability (or more precisely,
the ratio of permeability to dynamic viscosity) which is modelled as a
random field. Uncertainty in $a(\bsx,\omega)$ leads to uncertainty in
$q(\bsx,\omega)$ and $p(\bsx,\omega)$. Quantities of interest include for
example the breakthrough time of a plume of pollution moving through the
medium.

\subsubsection*{QMC for PDEs with random coefficients}

There is a huge literature on treating these PDEs with random coefficients
using various methods, see e.g., the surveys \cite{CD15,GWZ14,SG11} and
the references therein. Here we are interested in the application of
\emph{quasi-Monte Carlo} (\emph{QMC}) \emph{methods}, which are
equal-weight quadrature rules for high dimensional integrals, see e.g.,
\cite{DKS13,DP10,Lem09,LP14,Nie92,Nuy14,SJ94}.

QMC methods are still relatively new to these PDE problems. It began with
the 2011 paper \cite{GKNSS11} which included comprehensive numerical
experiments showing promising QMC results, but without any theoretical
justification. The first fully justified theory was provided in the 2012
paper \cite{KSS12}, and this has lead to a flood of research activities.
We will follow the recent survey \cite{KN16} to provide a high-level
overview of how QMC theory can be applied to PDEs with random
coefficients. The survey \cite{KN16} covered the detailed analysis from
six papers \cite{DKLNS14,DKLS16,GKNSSS15,KSSSU,KSS12,KSS15} in a unified
view. Different algorithms have been analyzed: \emph{single-level} vs
\emph{multi-level}, \emph{deterministic} vs \emph{randomized}, and
\emph{first order} vs \emph{higher order}, and they were considered under
different models for the randomness as we explain below.

It is popular to assume that $a(\bsx,\omega)$ is a \emph{lognormal} random
field, that is, $\log (a(\bsx,\omega))$ is a Gaussian random field on the
spatial domain $D$ with a specified mean and covariance function. Then one
can use the \emph{Karhunen--Lo\`eve} (\emph{KL}) \emph{expansion} to write
$\log (a(\bsx,\omega))$ as an infinite series parametrised by a sequence
$y_j = y_j(\omega)$, $j\ge 1$, of i.i.d.\ standard normal random numbers
from $\bbR$. Aside from the lognormal case, often the simpler
\emph{uniform} case is considered, where $a(\bsx,\omega)$ is written as an
infinite series that depends linearly on a sequence $y_j = y_j(\omega)$,
$j\ge 1$, of i.i.d.\ uniform random numbers from a bounded interval of
$[-1,1]$ or $[-\frac{1}{2},\frac{1}{2}]$. In both the lognormal and
uniform cases the infinite series is truncated in practice to, say, $s$
terms. The expected value of any quantity of interest is then approximated
by an $s$-dimensional integral with respect to the parameters $y_j$, which
can in turn be approximated by QMC methods, combined with finite element
methods for solving the PDE.

The six papers surveyed in \cite{KN16} all followed this KL-based general
direction. With respect to the QMC method they can be either \emph{first
order} or \emph{higher order}, which refers to the rate of convergence
being close to $\calO(n^{-1})$ or $\calO(n^{-\alpha})$, $\alpha>1$, with
$n$ being the number of integrand evaluations. With respect to the
approximation of the integrand function they can be either
\emph{single-level} or \emph{multi-level}, which refers to how spatial
discretization and dimension truncation are performed. A summary of the
results is given in the table below:

\begin{center}
 \begin{tabular}{lc@{\quad}c}
 & Uniform case & Lognormal case \\
 \hline
  First order single-level analysis & \cite{KSS12} & \cite{GKNSSS15} \\
  First order multi-level analysis & \cite{KSS15} & \cite{KSSSU} \\
  Higher order single-level analysis & \cite{DKLNS14} & \\
  Higher order multi-level analysis & \cite{DKLS16} & \\
 \end{tabular}
\end{center}

\noindent The first order results \cite{KSS12,KSS15} and
\cite{GKNSSS15,KSSSU} are based on \emph{randomly shifted lattice rules}
and are accompanied by probabilistic error bounds. The higher order
results \cite{DKLNS14,DKLS16} are based on \emph{interlaced polynomial
lattice rules} and are accompanied by deterministic error bounds. The
lognormal results \cite{GKNSSS15,KSSSU} require a non-standard function
space setting for integrands with domain $\bbR^s$. A key feature in all
these analysis is that the QMC error bounds are independent of the number
of integration variables $s$. There is as yet no satisfactory QMC theory
that can give higher order convergence for the lognormal case with error
bound independent of~$s$.

\subsubsection*{Plan of this article}

In Section~\ref{sec:overview} we provide an overview of the different
settings and algorithms covered in the survey \cite{KN16}, with the goal
to convey the overall picture while keeping the exposition as simple and
accessible as possible. In Section~\ref{sec:tut} we take a change of pace
and style to give a step-by-step tutorial of the required analysis for the
uniform case with first order QMC rules. That is, we zoom in and focus on
the essence of the paper \cite{KSS12} in such a way that the tutorial can
be used to extend the analysis to other cases by interested readers. Then
in Section~\ref{sec:key} we zoom out again and continue to provide
insights to the key analysis required for the six papers surveyed in
\cite{KN16}. In Section~\ref{sec:soft} we briefly discuss the software
accompanying \cite{KN16}. Finally in Section~\ref{sec:conc} we give a
short conclusion.

\subsubsection*{Beyond the survey}

There have been many developments beyond the scope of the survey
\cite{KN16}.

Instead of using the KL expansion, in the lognormal case one can sample
the random field only at a discrete set of points with respect to the
covariance matrix inherited from the given covariance function of the
continuous field. The random field is then represented exactly at these
points, thus eliminating completely the truncation error associated with
the KL-based approach. (Note that interpolation may be required at the
finite element quadrature nodes.) The resulting large matrix factorization
problem could potentially be handled by \emph{circulant embedding} and
FFT, if the covariance function is \emph{stationary} and the grid is
regular, see \cite{DN97}. In fact, this was the approach taken in the
first QMC paper for PDEs with random coefficients \cite{GKNSS11}, and the
corresponding analysis is being considered in
\cite{GKNSS-paper1,GKNSS-paper2}.

Another way to tackle the large matrix factorization is to make use of
\emph{$H$-matrix techniques}, see \cite{Hac15}, and this has been
considered in \cite{FKS}.

The uniform framework can be extended from the elliptic PDE to the general
framework of \emph{affine} parametric operator equations, see \cite{Sch12}
as well as \cite{DKLNS14,DKLS16}. A different QMC theory for the lognormal
case is offered in \cite{HPS17}. Further PDE computations with higher
order QMC are reported in \cite{GanSch16}, and with multi-level and
multi-index QMC in \cite{RNV}. QMC has also been applied to PDEs on the
sphere \cite{Leg12}, holomorphic equations \cite{DLS16}, Bayesian
inversion \cite{DGLS,SST17}, stochastic wave propagation
\cite{GH15,GanKS}, and eigenproblems \cite{GGKSS}.

Moreover, there has been some significant development in the use of
functions with local support in the expansions of $a(\bsx,\omega)$ which
leads to a simplified norm estimate for the integrand and a reduced
construction cost (pre-computation) for QMC, see \cite{GHS,HerSch,Kaz}.

\section{Overview} \label{sec:overview}

Throughout this article we refer to the number of integration variables
$s$ as the \emph{stochastic dimension}, which can be in the hundreds or
thousands or more (and controls the truncation error), in contrast to the
\emph{spatial dimension} $d$ which is just $1$, $2$ or $3$.

\subsection{Uniform vs Lognormal}

For a given parameter $\bsy$ we consider the parametric elliptic Dirichlet
problem
\begin{equation}\label{eq:strong}
 - \nabla \cdot (a(\bsx,\bsy)\,\nabla u(\bsx,\bsy))
 \,=\, \kappa(\bsx) \quad\mbox{for $\bsx$ in $D$}\,,
 \quad
 u(\bsx,\bsy) \,=\, 0 \quad\mbox{for $\bsx$ on $\partial D$}\,,
\end{equation}
for domain $D\subset\bbR^d$ a bounded, convex, Lipschitz polyhedron with
boundary $\partial D$, where the spatial dimension $d=1,2$, or $3$ is
assumed given and fixed. The differential operators in~\eqref{eq:strong}
are understood to be with respect to the physical variable $\bsx$ which
belongs to $D$. The parametric variable $\bsy = (y_j)_{j\geq 1}$ belongs
to either a bounded or unbounded domain, depending on which of the two
popular formulations of the parametric coefficient $a(\bsx,\bsy)$ is being
considered.

\subsubsection*{Uniform case}

In the uniform case, we assume that the $y_j$ are independent and
uniformly distributed on $[-\tfrac{1}{2},\tfrac{1}{2}]$, and
\begin{align} \label{eq:axy-unif}
  a(\bsx,\bsy)
  \,=\, a_0(\bsx) + \sum_{j\geq 1} y_j\, \psi_j(\bsx)\,,
  %\qquad
  %y_j \sim \calU[-\tfrac{1}{2},\tfrac{1}{2}]\,,
\end{align}
with $0 < a_{\min}\le a(\bsx,\bsy)\le a_{\max}<\infty$ for all $\bsx$ and
$\bsy$. We need further assumptions on $a_0$ and $\psi_j$, see \cite{KN16}
for details. Here we mention only one important assumption that there
exists $p_0\in (0,1)$ such that
\begin{equation} \label{eq:sum-p0}
 \sum_{j\ge 1}  \|\psi_j\|^{p_0}_{L_\infty} \,<\, \infty\,.
\end{equation}
The value of $p_0$ reflects the rate of decay of the fluctuations in
\eqref{eq:axy-unif}; later we will see that it directly affects the QMC
convergence rate.

Our goal is to compute the integral, i.e., the expected value, with
respect to $\bsy$, of a bounded linear functional $G$ applied to the
solution $u(\cdot,\bsy)$ of the PDE \eqref{eq:strong}
\begin{align} \label{eq:int-unif}
  \int_{[-\tfrac{1}{2},\tfrac{1}{2}]^\bbN} G(u(\cdot,\bsy))\,\rd\bsy
  \,:=\, \lim_{s\to\infty} \int_{[-\tfrac{1}{2},\tfrac{1}{2}]^s}
        G(u(\cdot,(y_1,\ldots,y_s,0,0,\ldots)))\,\rd y_1\cdots\rd y_s\,.
\end{align}

\subsubsection*{Lognormal case}

In the lognormal case, we assume that the $y_j$ are independent standard
normal random numbers on $\bbR$, and
\begin{align} \label{eq:axy-logn}
  a(\bsx,\bsy) \,=\,
  a_0(\bsx)\exp\bigg(\sum_{j\ge 1}
  y_j\,\sqrt{\mu_j}\,\xi_j(\bsx)\bigg)\,,
  %\qquad
  %y_j \sim \calN(0,1)\,,
\end{align}
where $a_0(\bsx)>0$, the $\mu_j>0$ are non-increasing, and the $\xi_j$ are
orthonormal in $L_2(D)$. This can arise from the KL expansion in the case
where $\log(a)$ is a stationary Gaussian random field with a specified
mean and covariance function; a popular choice is the \emph{Mat\'ern}
covariance.

Our goal now is the integral of $G(u(\cdot,\bsy))$ over $\bsy\in
\bbR^\bbN$ with a countable product Gaussian measure $\mu_G(\rd\bsy)$
(formally, we restrict the domain to some $Y\subset\bbR^\bbN$ with full
measure $\mu_G(Y) = 1$, but we omit this in the notation)
\begin{align} \label{eq:int-logn}
 \int_{\bbR^\bbN} G(u(\cdot,\bsy))\,\prod_{j\ge 1} \phi_{\rm nor}(y_j)\,\rd\bsy
 \,=\, \int_{[0,1]^\bbN} G(u(\cdot,\Phi^{\mbox{-}1}_{\rm nor}(\bsw)))\,\rd\bsw\,,
\end{align}
where $\phi_{\rm nor}(y) := \exp(-y^2/2)/\sqrt{2\pi}$ is the univariate
standard normal probability density function, while $\Phi^{\mbox{-}1}_{\rm
nor}$ denotes the inverse of the corresponding cumulative distribution
function, and is applied component-wise to a vector. The transformed
integral over the unit cube on the right-hand side of~\eqref{eq:int-logn}
is obtained by the change of variables $\bsy = \Phi^{\mbox{-}1}_{\rm
nor}(\bsw)$.

\subsection{Single-level vs Multi-level}

\subsubsection*{Single-level algorithms}

We approximate the integral~\eqref{eq:int-unif} or~\eqref{eq:int-logn} in
three steps:
\begin{enumerate}
\item [i.] Dimension truncation: the infinite sum
    in~\eqref{eq:axy-unif} or~\eqref{eq:axy-logn} is truncated to $s$
    terms.
\item [ii.] Finite element discretization: the PDE~\eqref{eq:strong}
    in weak formulation (see \eqref{eq:weak} below) is solved using a
    finite element method with meshwidth $h$.
\item [iii.] QMC quadrature: the integral of the finite element
    solution for the truncated problem is estimated using a
    deterministic or randomized QMC method.
\end{enumerate}

The deterministic version of this algorithm is
\begin{align} \label{eq:SL-DET}
  \frac{1}{n} \sum_{i=1}^{n} G(u^s_h(\cdot,\bsy_i))\,,
  \qquad
  \bsy_i
  \,=\,
  \begin{cases}
  \bst_i - \tfrac{\boldsymbol{1}}{\boldsymbol{2}} & \mbox{for uniform}, \\
  \Phi^{\mbox{-}1}_{\rm nor}(\bst_i) & \mbox{for lognormal},
  \end{cases}
\end{align}
where $\bst_1,\ldots,\bst_n \in [0,1]^s$ are $n$ QMC points from the
$s$-dimensional standard unit cube. In the uniform case, these points are
translated to the unit cube $[-\tfrac{1}{2},\frac{1}{2}]^s$. In the
lognormal case, these points are mapped to the Euclidean space $\bbR^s$ by
applying the inverse of the cumulative normal distribution function
component-wise.

A randomized version of this algorithm with \emph{random shifting} is
given by
\begin{align} \label{eq:SL-RAN}
  \frac{1}{r} \sum_{k=1}^{r}
  \frac{1}{n} \sum_{i=1}^{n} G(u^s_h(\cdot,\bsy_{i,k}))\,,
  \qquad
  \bsy_{i,k}
  \,=\,
  \begin{cases}
  \{\bst_i + \bsDelta_k\} - \tfrac{\boldsymbol{1}}{\boldsymbol{2}} & \mbox{for uniform}, \\
  \Phi^{\mbox{-}1}_{\rm nor}(\{\bst_i + \bsDelta_k\}) & \mbox{for lognormal},
  \end{cases}
\end{align}
where $\bst_1,\ldots,\bst_n \in [0,1]^s$ are $n$ QMC points as above, and
$\bsDelta_1,\ldots,\bsDelta_r \in [0,1]^s$ are $r$ independent
\emph{random shifts} generated from the uniform distribution on $[0,1]^s$.
The braces in $\{\bst_i + \bsDelta_k\}$ mean that we take the fractional
part of each component in the vector $\bst_i + \bsDelta_k$. Other
randomization strategies can be used analogously but need to be chosen
appropriately to preserve the special properties of the QMC points.
Randomized algorithms have the advantages of being unbiased as well as
providing a practical error estimate.

\subsubsection*{Multi-level algorithms}

The general concept of \emph{multi-level} can be explained as follows (see
e.g., \cite{Gil15}): if we denote the integral~\eqref{eq:int-unif}
or~\eqref{eq:int-logn} by $I_\infty$ and define a sequence $I_0, I_1,
\ldots$ of approximations converging to $I_\infty$, then we can write
$I_\infty$ as a telescoping sum $I_\infty = (I_\infty - I_L) +
\sum_{\ell=0}^L (I_\ell - I_{\ell-1})$, $I_{-1}:=0$, and then apply
different quadrature rules to the differences $I_\ell - I_{\ell-1}$, which
we anticipate to get smaller as $\ell$ increases. Here we define $I_\ell$
to be the integral of $G(u^{s_\ell}_{h_\ell})$ corresponding to the finite
element solution with meshwidth $h_\ell$, for the truncated problem with
$s_\ell$ terms, where $1 \le s_0 \le s_1 \le s_2 \le \cdots \le s_L \le
\cdots$ and $h_0 \ge h_1 \ge h_2 \ge \cdots \ge h_L \ge \cdots > 0$, so
that $I_\ell$ becomes a better approximation to $I_\infty$ as $\ell$
increases.

The deterministic version of our multi-level algorithm takes the form
(remembering the linearity of $G$)
\begin{align} \label{eq:ML-DET}
  \sum_{\ell=0}^L \bigg(\frac{1}{n_\ell} \sum_{i=1}^{n_\ell}
                 G((u^{s_\ell}_{h_\ell}-u^{s_{\ell-1}}_{h_{\ell-1}})
                 (\cdot,\bsy_i^\ell))\bigg)\,,
  \qquad
  \bsy_i^\ell \,=\,
  \begin{cases}
  \bst_i^\ell - \tfrac{\boldsymbol{1}}{\boldsymbol{2}} & \mbox{for uniform}, \\
  \Phi^{\mbox{-}1}_{\rm nor}(\bst_i^\ell) & \mbox{for lognormal},
  \end{cases}
\end{align}
where we apply an $s_\ell$-dimensional QMC rule with $n_\ell$ points
$\bst_1^\ell,\ldots,\bst_{n_\ell}^\ell \in [0,1]^{s_\ell}$ to the
integrand $G(u^{s_\ell}_{h_\ell}-u^{s_{\ell-1}}_{h_{\ell-1}})$, and we
define $u^{s_{-1}}_{h_{-1}} :=0$.

The corresponding randomized version can be obtained analogously to
\eqref{eq:SL-RAN} by taking $r_\ell$ random shifts at each level, noting
that all shifts from all levels should be independent.

\subsection{First-order vs Higher-order}

Up to this point we have said very little about QMC methods, other than
noting that they are equal-weight quadrature rules as seen in
\eqref{eq:SL-DET}. Actually, we will not say much about QMC methods in
this article at all. In this subsection we will mention three different
QMC theoretical settings which have been used for PDEs applications,
giving just enough details in the first setting needed for the tutorial in
Section~\ref{sec:tut}. These three settings are discussed in slightly more
detail in \cite{KN-plenary} in this volume, and more comprehensively in
\cite{KN16}; see also the references in these papers.

\subsubsection*{First order QMC over the unit cube -- randomly shifted lattice rules for weighted Sobolev spaces}

Suppose we wish to approximate the $s$-dimensional integral over the unit
cube $[0,1]^s$
\begin{align} \label{eq:int1}
 \int_{[0,1]^s} f(\bsy)\, \rd\bsy\,,
\end{align}
where the integrand $f$ belongs to a \emph{weighted Sobolev space of
smoothness one}, with the \emph{unanchored} norm defined by (see e.g.,
\cite{SWW04})
\begin{align} \label{eq:norm1}
  \|f\|_\bsgamma
  \,=\,
  \Bigg[
  \sum_{\setu\subseteq\{1:s\}}
  \frac{1}{\gamma_\setu}
  \int_{[0,1]^{|\setu|}}
  \bigg(\int_{[0,1]^{s-|\setu|}}
  \frac{\partial^{|\setu|}f}{\partial \bsy_\setu}(\bsy)
  \,\rd\bsy_{\{1:s\}\setminus\setu}
  \bigg)^2
  \rd\bsy_\setu
  \Bigg]^{1/2}.
\end{align}
Here $\{1:s\}$ is a shorthand notation for the set of indices
$\{1,2,\ldots,s\}$, $(\partial^{|\setu|}f)/(\partial \bsy_\setu)$ denotes
the mixed first derivative of $f$ with respect to the ``active'' variables
$\bsy_\setu = (y_j)_{j\in\setu}$, while $\bsy_{\{1:s\}\setminus\setu} =
(y_j)_{j\in\{1:s\}\setminus\setu}$ denotes the ``inactive'' variables.
There is a weight parameter $\gamma_\setu\ge 0$ associated with each
subset of variables~$\bsy_\setu$ to moderate the relative importance
between the different sets of variables. We denote the weights
collectively by $\bsgamma$.

In this setting we pair the weighted Sobolev space with \emph{randomly
shifted lattice rules}; the complete theory can be found in \cite{DKS13}.
They approximate the integral \eqref{eq:int1} by
\[
  \frac{1}{n} \sum_{i=1}^n f(\bst_i),
  \qquad \bst_i \,=\, \left\{\frac{i\,\bsz}{n} + \bsDelta\right\}\,,
\]
where $\bsz\in \bbZ^s$ is known as the \emph{generating vector},
$\bsDelta$ is a \emph{random shift} drawn from the uniform distribution
over $[0,1]^s$, and as in \eqref{eq:SL-RAN} the braces indicate that we
take the fractional parts of each component in a vector. It is known that
good generating vectors can be obtained using a \emph{CBC construction}
(\emph{component-by-component construction}), determining the components
of $\bsz$ one at a time sequentially, to achieve first order convergence
in this setting, where the implied constant can be independent of $s$
under appropriate conditions on the weights $\bsgamma$.

Specifically, if $n$ is a power of $2$ then we know that the CBC
construction yields the root-mean-square error bound (with respect to the
uniform random shift), for all $\lambda\in (1/2,1]$,
\begin{align} \label{eq:error1}
  \mbox{r.m.s. error}
  &\,\le\,
  \Bigg(
  \frac{2}{n}
  \sum_{\emptyset\ne\setu\subseteq\{1:s\}} \gamma_\setu^\lambda\,
  [\vartheta(\lambda)]^{|\setu|}
  \Bigg)^{1/(2\lambda)}
  \,\|f\|_\bsgamma\,,
\end{align}
where $\vartheta(\lambda) := 2\zeta(2\lambda)/(2\pi^2)^\lambda$, with
$\zeta(a) := \sum_{k=1}^\infty k^{-a}$ denoting the Riemann zeta function.
A similar result holds for general $n$. The best rate of convergence
clearly comes from choosing $\lambda$ close to $1/2$.

We need some structure in the weights $\bsgamma$ for the CBC construction
cost to be feasible in practice. \emph{Fast} CBC algorithms (using FFT)
can find a generating vector of a good $n$-point lattice rule in $s$
dimensions in $\calO(s\,n\log n)$ operations in the case of \emph{product
weights}, and in $\calO(s\,n\log n + s^2\,n)$ operations in the case of
\emph{POD weights} (see \eqref{eq:step13} ahead).

\subsubsection*{First order QMC over $\bbR^s$}

We can pair randomly shifted lattice rules with a special function space
setting over $\bbR^s$ to achieve first order convergence. The norm in this
function space setting includes some additional weight functions to
control the behavior of the derivatives of the functions as the components
go to $\pm \infty$. The root-mean-square error bound takes the same form
as \eqref{eq:error1}, but with a different definition of the norm and
$\vartheta(\lambda)$.

\subsubsection*{Higher order QMC over the unit cube}

We can pair a family of QMC methods called \emph{interlaced polynomial
lattice rules} with another special function space setting over the unit
cube to achieve higher order convergence. The norm in this function space
setting involves higher order mixed derivatives of the functions. The key
advantage of this family of QMC methods over other higher order QMC
methods is that, in the cost of finding a generating vector which achieves
the best theoretical convergence rate, the \emph{order} or the
\emph{interlacing factor} appears as a multiplying factor rather than
sitting in the exponent of the number of points $n$.

\section{Tutorial} \label{sec:tut}

We conclude from the error bound \eqref{eq:error1} that the first step in
applying QMC theory is to estimate the norm of the practical integrand. We
see from \eqref{eq:SL-DET}, \eqref{eq:SL-RAN}, and \eqref{eq:ML-DET} that
this means we need to estimate the norms
\[
  \|G(u^s_h)\|_{\bsgamma}
  \qquad\mbox{and}\qquad
  \|G(u^{s_\ell}_{h_\ell}-u^{s_{\ell-1}}_{h_{\ell-1}})\|_{\bsgamma}\,,
\]
for the single-level and the multi-level algorithms, respectively.

In this section we provide a step-by-step tutorial on the analysis for the
single-level algorithm in the uniform case with first order QMC rules.

\subsubsection*{Differentiate the PDE}

\begin{enumerate}
\item %1
We start with the variational formulation of the PDE
\eqref{eq:strong}: find $u(\cdot,\bsy)\in H^1_0(D)$ such that
\begin{equation} \label{eq:weak}
 \int_D a(\bsx,\bsy)\,\nabla u(\bsx,\bsy)\cdot\nabla w(\bsx)\,\rd\bsx
 \,=\, \int_D \kappa(\bsx)\,w(\bsx)\,\rd\bsx \quad\forall w\in H^1_0(D)\,.
\end{equation}
Here we consider the Sobolev space $H^1_0(D)$ of functions which
vanish on the boundary of $D$, with norm $\|w\|_{H^1_0} := \|\nabla
w\|_{L_2}$, and together with the dual space $H^{\mbox{-}1}(D)$ and
pivot space $L_2(D)$.
\smallskip

\item %2
We take the mixed partial derivatives $\partial^{\bsnu}$ with respect
to $\bsy$ with multi-index $\bsnu\ne\bszero$ (i.e., we differentiate
$\nu_j$ times with respect to $y_j$ for each $j$) on both sides of
\eqref{eq:weak} to obtain
\begin{equation} \label{eq:step2}
 \int_D \partial^{\bsnu}\Big(
 a(\bsx,\bsy)\, \nabla u(\bsx,\bsy)\cdot\nabla
 w(\bsx)\Big)\,\rd\bsx \,=\, 0 \quad\forall w\in H^1_0(D)\,.
\end{equation}
We can move the derivatives inside the integrals because they operate
on different variables $\bsy$ and $\bsx$, respectively. The right-hand
side vanishes because it does not depend on $\bsy$. \smallskip

\item %3
Next we apply the Leibniz product rule on the left-hand side of
\eqref{eq:step2} to obtain
\begin{equation} \label{eq:step3}
 \int_D \bigg( \sum_{\bsm\le\bsnu}
 \binom{\bsnu}{\bsm}(\partial^{\bsm}a)(\bsx,\bsy)\,
 \nabla (\partial^{\bsnu-\bsm} u)(\bsx,\bsy)\cdot\nabla
 w(\bsx)\bigg)\,\rd\bsx \,=\, 0 \quad\forall w\in H^1_0(D)\,,
\end{equation}
where the sum is over all multi-indices $\bsm$ satisfying
$\bsm\le\bsnu$ (i.e., $m_j\le\nu_j$ for all $j$), and
$\binom{\bsnu}{\bsm} := \prod_{j\ge 1} \binom{\nu_j}{m_j}$. So far we
have made no use of any assumption on $a(\bsx,\bsy)$. \smallskip

\item %4
For the uniform case, it is easy to see from the formula
\eqref{eq:axy-unif} of $a(\bsx,\bsy)$ that
\begin{equation} \label{eq:step4}
  (\partial^{\bsm}a)(\bsx,\bsy) \,=\,
 \begin{cases}
 a(\bsx,\bsy) &\mbox{if } \bsm=\bszero\,, \\
 \psi_j(\bsx) &\mbox{if } \bsm=\bse_j\,, \\
 0 &\mbox{otherwise}\,,
 \end{cases}
\end{equation}
where $\bse_j$ denotes the multi-index whose $j$th component is $1$
and all other components are $0$. Essentially, due to the linearity of
$a$ with respect to each $y_j$, if we differentiate once then we
obtain $\psi_j$, and if we differentiate a second time with respect to
any variable we get $0$. \smallskip

\item %5
Substituting \eqref{eq:step4} into \eqref{eq:step3} and separating out
the $\bsm=\bszero$ term, we obtain
\begin{align} \label{eq:step5}
 &\int_D
    a(\bsx,\bsy)\,\nabla (\partial^{\bsnu}
    u)(\bsx,\bsy)\cdot \nabla w(\bsx)\,\rd\bsx \nonumber\\
 &\,=\, - \sum_{j\ge 1} \nu_j \int_D
   \psi_j(\bsx)\, \nabla (\partial^{\bsnu-\bse_j}
    u)(\bsx,\bsy) \cdot\nabla w(\bsx)\,\rd\bsx\quad\forall
    w\in H^1_0(D)\,.
\end{align}

\item %6
Note that \eqref{eq:step5} holds for all test functions in $H^1_0(D)$.
We now take the particular choice of $w = (\partial^{\bsnu}
u)(\cdot,\bsy)$ (yes, it is allowed to depend on $\bsy$) in
\eqref{eq:step5}. Applying $a(\bsx,\bsy)\ge a_{\min}$ to the left-hand
side, and $|\psi_j(\bsx)|\le \|\psi_j\|_{L_\infty}$ and the
Cauchy-Schwarz inequality to the right-hand side, we obtain
\begin{align} \label{eq:step6}
 &a_{\min}
 \int_D |\nabla (\partial^{\bsnu} u)(\bsx,\bsy)|^2\,\rd\bsx \\
 &\,\le\,
 \sum_{j\ge 1} \nu_j\, \|\psi_j\|_{L_\infty} \bigg(\int_D
 |\nabla (\partial^{\bsnu-\bse_j}
 u)(\bsx,\bsy)|^2\,\rd\bsx\bigg)^{1/2} \bigg(\int_D |\nabla
 (\partial^{\bsnu}
 u)(\bsx,\bsy)|^2\,\rd\bsx\bigg)^{1/2}\,. \nonumber
\end{align}

\item %7
Canceling one common factor from both sides of \eqref{eq:step6} and
then dividing through by $a_{\min}$, we obtain the recurrence
\begin{align} \label{eq:step7}
\|\nabla
    (\partial^{\bsnu} u)(\cdot,\bsy)\|_{L_2} \,\le\, \sum_{j\ge
    1} \nu_j\, b_j\, \|\nabla (\partial^{\bsnu-\bse_j}
    u)(\cdot,\bsy)\|_{L_2}\,, \qquad
    b_j \,:=\,
    \frac{\|\psi_j\|_{L_\infty}}{a_{\min}}\,.
\end{align}

\ifspringer
\clearpage %%%%%%%%%%%%%%%%%%%%%%%%%%%%%%%%%%%%%%%%%%%%%%%%%%%%%%%%% FIX ME
\fi

\item %8
Finally we prove by induction that
\begin{align} \label{eq:step8}
 \|\nabla (\partial^{\bsnu} u)(\cdot,\bsy)\|_{L_2}
 \,\le\,
 |\bsnu|!\,\bsb^{\bsnu}\,
 \frac{\|\kappa\|_{H^{\mbox{-}1}}}{a_{\min}}\,,
\end{align}
where $|\bsnu| := \sum_{j\ge 1} \nu_j$ and $\bsb^\bsnu := \prod_{j\ge
1} b_j^{\nu_j}$.

\begin{enumerate}
\item \textsc{Base step.} %
We return to the variational form \eqref{eq:weak} and take $w =
u(\cdot,\bsy)$. Applying $a(\bsx,\bsy)\ge a_{\min}$ to the
left-hand side and estimating the right-hand side using duality
pairing $|\langle \kappa, u(\cdot,\bsy)\rangle| \le
\|\kappa\|_{H^{\mbox{-}1}}\,\|u(\cdot,\bsy)\|_{H^1_0}$, we obtain
\[
  a_{\min}\,\|\nabla u(\cdot,\bsy)\|_{L_2}^2
  \,\le\, \|\kappa\|_{H^{\mbox{-}1}}\,\|\nabla u(\cdot,\bsy)\|_{L_2}\,,
\]
which can be rearranged to yield the case $\bsnu = \bszero$ in
\eqref{eq:step8}.
\smallskip

\item \textsc{Induction step.} %
As the induction hypothesis, we assume that \eqref{eq:step8} holds
for all multi-indices of order $< |\bsnu|$. Then we have
\[
 \|\nabla (\partial^{\bsnu-\bse_j}u)(\cdot,\bsy)\|_{L_2}
 \,\le\, |\bsnu-\bse_j|!\,\bsb^{\bsnu-\bse_j}\,
 \frac{\|\kappa\|_{H^{\mbox{-}1}}}{a_{\min}}\,.
\]
Substituting this into \eqref{eq:step7} and noting that
$\nu_j\,|\bsnu-\bse_j|! = |\bsnu|!$ and $b_j\,\bsb^{\bsnu-\bse_j}
= \bsb^\bsnu$, we obtain \eqref{eq:step8} and conclude the
induction.
\end{enumerate}
\end{enumerate}

\subsubsection*{Estimate the norm}

\begin{enumerate}
\setcounter{enumi}{8}
\item %9
We want to estimate the norm $\|G(u^s_h)\|_\bsgamma$. We see from the
definition of the norm in \eqref{eq:norm1} that we need to obtain
estimates on the mixed first derivatives of $G(u^s_h(\cdot,\bsy))$
with respect to $\bsy$. Using linearity and boundedness of $G$, we
have
\begin{align} \label{eq:step9}
 \bigg|\frac{\partial^{|\setu|}}{\partial
    \bsy_\setu}
   G(u^{s}_h(\cdot,\bsy))\bigg|
   \,=\,
   \bigg|G\bigg(\frac{\partial^{|\setu|}}{\partial \bsy_\setu}
   u^{s}_h(\cdot,\bsy)\bigg)\bigg| \le
   \|G\|_{H^{\mbox{-}1}}\bigg\|\frac{\partial^{|\setu|}}{\partial \bsy_\setu}
   u^{s}_h(\cdot,\bsy)\bigg\|_{H^1_0}\,.
\end{align}

\item %10
We can repeat the above proof of \eqref{eq:step8} for the truncated
finite element solution $u^s_h$ instead of the true solution $u$. Then
we restrict the result to mixed first derivatives (i.e., $\nu_j\le 1$
for all $j$) and deduce that
\begin{align} \label{eq:step10}
 \bigg\|\frac{\partial^{|\setu|}}{\partial
 \bsy_\setu}
   u^{s}_h(\cdot,\bsy)\bigg\|_{H^1_0} \,\le\, |\setu|!
   \bigg(\prod_{j\in\setu} b_j\bigg)
   \frac{\|\kappa\|_{H^{\mbox{-}1}}}{a_{\min}}, \quad \setu\subseteq\{1:s\}\,.
\end{align}

\item %11
Combining \eqref{eq:step9} with \eqref{eq:step10} and substituting the
upper bound into the definition of the norm \eqref{eq:norm1}, we
conclude that
\begin{align} \label{eq:step11}
 \|
 G(u^s_h)\|_{\bsgamma}
      \,\le\, \frac{\|\kappa\|_{H^{\mbox{-}1}}\|G\|_{H^{\mbox{-}1}}}{a_{\min}}
  \Bigg(\sum_{\setu\subseteq\{1:s\}} \frac{(|\setu|!)^2
  \prod_{j\in\setu} b_j^2}{\gamma_\setu} \Bigg)^{1/2}\,.
\end{align}
\end{enumerate}

\subsubsection*{Choose the weights}

\begin{enumerate}

\setcounter{enumi}{11}
\item %12
Now we apply the upper bound on the norm \eqref{eq:step11} in the
error bound for randomly shifted lattice rules \eqref{eq:error1}, to
yield (leaving out some constants as indicated by $\lesssim$) for all
$\lambda\in (1/2,1]$,
\begin{align} \label{eq:step12}
 \mbox{r.m.s.\ error} \,\lesssim\, \Bigg( \frac{2}{n}
  \sum_{\setu\subseteq\{1:s\}} \gamma_\setu^\lambda\,
  [\vartheta(\lambda)]^{|\setu|} \Bigg)^{1/(2\lambda)} \Bigg(
  \sum_{\setu\subseteq\{1:s\}} \frac{(|\setu|!)^2
  \prod_{j\in\setu} b_j^2}{\gamma_\setu} \Bigg)^{1/2}\,.
\end{align}

\item %13
With elementary calculus, for any $\lambda$, we can minimize the the
upper bound in \eqref{eq:step12} with respect to the weights
$\gamma_\setu$ to yield the formula
\begin{align} \label{eq:step13}
 \gamma_\setu \,=\, \bigg(|\setu|!\prod_{j\in\setu}
  \frac{b_j}{\sqrt{\vartheta}(\lambda)} \bigg)^{2/(1+\lambda)}.
\end{align}
This form of weights is called \emph{product and order dependent
weights}, or \emph{POD weights} in short, because of the presence of
some product factors as well as the cardinality of $\setu$. \smallskip

\item %14
We substitute \eqref{eq:step13} into \eqref{eq:step12} and simplify
the expression to
\begin{align} \label{eq:step14}
 &\mbox{r.m.s.\ error} \\
 &\quad\lesssim\, \Bigg( \frac{2}{n}
  \Bigg)^{1/(2\lambda)} \Bigg[
  \sum_{\setu\subseteq\{1:s\}} \Bigg( |\setu|!
  \prod_{j\in\setu} \Big(b_j\, [\vartheta(\lambda)]^{1/(2\lambda)}\Big)\Bigg)^{2\lambda/(1+\lambda)} \
  \Bigg]^{(1+\lambda)/(2\lambda)}\,. \nonumber
\end{align}

\item %15
We now derive a condition on $\lambda$ for which the sum in
\eqref{eq:step14} is bounded independently of $s$. In an abstract
form, we have
\[
 \sum_{\setu\subseteq\{1:s\}} \bigg(|\setu|! \prod_{j\in\setu} \alpha_j\bigg)^k
 \,=\, \sum_{\ell=0}^s (\ell!)^k
 \!\!
 \sum_{\setu\subseteq\{1:s\},\,|\setu|=\ell}
 \;\;
 \prod_{j\in\setu} \alpha_j^k
 \,\le\, \sum_{\ell=0}^s (\ell!)^{k-1} \bigg(\sum_{j=1}^s \alpha_j^k\bigg)^\ell\,,
\]
where the inequality holds because each term $\prod_{j\in\setu}
\alpha_j^k$ from the left-hand side of the inequality appears in the
expansion $(\sum_{j=1}^s \alpha_j^k)^\ell$ exactly $\ell!$ times and
yet the expansion contains other terms. The right-hand side is bounded
independently of~$s$ if $\sum_{j=1}^\infty \alpha_j^k<\infty$ and
$k<1$, which can be verified by the ratio test. In our case, we have
$k = 2\lambda/(1+\lambda)$ and $\sum_{j=1}^\infty \alpha_j^k =
[\vartheta(\lambda)]^{1/(1+\lambda)}\,\sum_{j=1}^\infty b_j^k <
\infty$ if $k\ge p_0$, where we recall that $b_j$ is defined in
\eqref{eq:step7} and $p_0$ is defined in \eqref{eq:sum-p0}. Hence we
require
\begin{align} \label{eq:step15}
  p_0 \,\le\, \frac{2\lambda}{1+\lambda} \,<\, 1
  \quad\Longleftrightarrow
  \quad
  \frac{p_0}{2-p_0} \,\le\, \lambda \,<\, 1\,.
\end{align}

\item %16
Clearly the best rate of convergence is obtained by taking $\lambda$
as small as possible. Combining the original constraint of $\lambda\in
(1/2,1]$ with \eqref{eq:step15}, we now take
\begin{align} \label{eq:step16}
 \lambda \,=\,
  \begin{cases}
  \displaystyle\frac{1}{2-2\delta} \mbox{ for } \delta\in (0,\tfrac{1}{2})
  & \mbox{when } p_0\in (0,\tfrac{2}{3}], \\[4mm]
  \displaystyle\frac{p_0}{2-p_0} & \mbox{when } p_0\in (\tfrac{2}{3},1).
  \end{cases}
\end{align}
\end{enumerate}

\subsubsection*{Fast CBC construction}

\begin{enumerate}
\setcounter{enumi}{16}
\item %17
The chosen weights \eqref{eq:step13} with $\lambda$ given by
\eqref{eq:step16} are then fed into the CBC construction to produce
tailored randomly shifted lattice rules that achieve a root-mean-square
error of order
\[
  n^{-\min(1/p_0-1/2,1-\delta)}, \quad \delta\in (0,\tfrac{1}{2})\,,
\]
with the implied constant independent of $s$, where $p_0$ is given by
\eqref{eq:sum-p0}. The fast CBC construction with POD weights can then
find a good generating vector in $\calO(s\,n\log n + s^2\,n)$
operations.
\end{enumerate}

\section{Key Analysis} \label{sec:key}

Having completed our embedded tutorial in the previous section, we now
continue to provide our overview of the analysis required in applying QMC
to PDEs with random coefficients.

\subsubsection*{Some hints at the technical difficulties for the multi-level analysis}

We have seen in the uniform case with the single-level algorithm that the
key is to estimate $\| G(u^s_h)\|_\bsgamma$, and this is achieved by
estimating (see \eqref{eq:step8} and \cite[Lemma~6.1]{KN16})
\[
  \| \nabla \partial^{\bsnu} u(\cdot,\bsy)\|_{L_2}.
\]
For the multi-level algorithm, the key estimate is $\| G(u^{
s_\ell}_{h_\ell}-u^{s_{\ell-1}}_{h_{\ell-1}})\|_{\bsgamma}$, and we need
to estimate in turn (see \cite[Lemmas~6.2--6.4]{KN16})
\[
 \| \Delta \partial^{\bsnu} u(\cdot,\bsy)\|_{L_2},
 \quad
 \| \nabla \partial^{\bsnu} (u-u_h)(\cdot,\bsy)\|_{L_2},
 \quad\mbox{and}\quad
  | \partial^{\bsnu} G((u-u_h)(\cdot,\bsy))|.
\]
All three bounds involve factors of the form $(|\bsnu|\!+\!a_1)!\,
\overline{\bsb}^{\bsnu}$ for $a_1\ge 0$ and a sequence $\overline{b}_j$
similar to the previously defined $b_j$. Assuming that both the forcing
term $\kappa$ and the linear functional $G$ are in $L_2(D)$, we obtain
that the second bound is of order $h$ and the third bound is of order
$h^2$. The difficulty is that we need to establish these regularity
estimates simultaneously in $\bsx$ and $\bsy$. We also use duality tricks
to gain on the convergence rate due to the linear functional $G$.

\subsubsection*{Some hints at the technical difficulties for the lognormal case}

For the lognormal case the argument is quite technical due to the more
complicated form of $a(\bsx,\bsy)$. In the single-level algorithm we need
to estimate (see \cite[Lemma~6.5]{KN16})
\[
\| \nabla \partial^{\bsnu} u(\cdot,\bsy)\|_{L_2}
  \quad\mbox{by first estimating}\quad
\| a^{1/2}(\cdot,\bsy)\,\nabla \partial^{\bsnu} u(\cdot,\bsy)\|_{L_2}\,.
\]
In the multi-level algorithm we need to estimate (see
\cite[Lemma~6.6]{KN16})
\[
 \| \Delta \partial^{\bsnu} u(\cdot,\bsy)\|_{L_2}
  \quad\mbox{by first estimating}\quad
\| a^{-1/2}(\cdot,\bsy) \nabla\cdot(a(\cdot,\bsy)\,\nabla \partial^{\bsnu} u(\cdot,\bsy))\|_{L_2}\,,
\]
and then estimate in turn (see \cite[Lemmas~6.7--6.8]{KN16})
\[
  \|a^{1/2}(\cdot,\bsy)\nabla \partial^{\bsnu} (u\!-\!u_h)(\cdot,\bsy)\|_{L_2}
  \quad\mbox{and}\quad
  | \partial^{\bsnu} G((u-u_h)(\cdot,\bsy))|\,.
\]
All bounds involve factors of the form
$J(\bsy)\,(|\bsnu|\!+\!a_1)!\,\bsbeta^{\bsnu}$ for $a_1\ge 0$ and some
sequence $\beta_j$, where $J(\bsy)$ indicates some factor depending on
$\bsy$ which is not present in the uniform case. The proofs are by
induction, and the tricky part is knowing what multiplying factor of
$a(\cdot,\bsy)$ should be included in the recursion. The growth of
$J(\bsy)$ needs to be taken into account when estimating the norm.

\subsubsection*{Summary of results}

Now we summarize and compare the results from
\cite{KSS12,KSS15,DKLNS14,DKLS16} for the uniform case:
\begin{align*}
 & \mbox{First-order single-level \cite{KSS12}} \\
 & \quad s^{-2(1/p_0-1)} + h^{t+t'} + n^{-\min(1/p_0-1/2,1-\delta)} \quad\mbox{(r.m.s.)} \\
 & \mbox{First-order multi-level \cite{KSS15}} \\
 & \quad s_L^{-2(1/p_0-1)} + h_L^{t+t'}
 + \sum_{\ell=0}^L n_\ell^{-\min(1/p_1-1/2,1-\delta)}
                   \big(\theta_{\ell-1}\,s_{\ell-1}^{-(1/p_0-1/p_1)}
                   + h_{\ell-1}^{t+t'}\big) \quad\mbox{(r.m.s.)} \\
 & \mbox{Higher-order single-level \cite{DKLNS14}} \\
 & \quad s^{-2(1/p_0-1)} + h^{t+t'} + n^{-1/p_0} \\
 & \mbox{Higher-order multi-level \cite{DKLS16}} \\
 & \quad s_L^{-2(1/p_0-1)} + h_L^{t+t'}
 + \sum_{\ell=0}^L n_\ell^{-1/p_t}
 \big(\theta_{\ell-1}\,s_{\ell-1}^{-(1/p_0-1/p_t)} +  h_{\ell-1}^{t+t'} \big)
\end{align*}
For the first-order results, the ``r.m.s.'' in brackets indicates that the
error is in the root-mean-square sense since we use a randomized QMC
method. The higher-order results are deterministic. Without giving the
full details, we simply say that the results include general parameters
$t$ and $t'$ for the regularity of $\kappa$ and $G$, respectively. Recall
that $p_0$ corresponds to the summability of $\|\psi_j\|_{L_\infty}$, see
\eqref{eq:sum-p0}. Here $p_1$ corresponds essentially to the summability
of $\|\nabla\psi_j\|_{L_\infty}$, while $p_t$ corresponds analogously to
higher derivatives of $\psi_j$. For the multi-level results we include the
analysis for potentially taking different $s_\ell$ at each level:
$\theta_{\ell-1}$ is $0$ if $s_{\ell} = s_{\ell-1}$ and is $1$ otherwise.

In the single-level algorithms, the error is the sum of three terms. In
the multi-level algorithms, we see the multiplicative effect between the
finite element error and the QMC error. However, comparing $p_1$ and $p_t$
with $p_0$, we see that multi-level algorithms need stronger regularity in
$\bsx$ than single-level algorithms.

Going from first-order to higher-order results, we see that the cap of
$n^{-(1-\delta)}$ is removed. We also see a gain of an extra factor of
$n^{-1/2}$; this benefit appears to arise from the switch of function
space setting to a non-Hilbert space.

The error versus cost analysis depends crucially on the cost assumptions.
For the single-level algorithms, we simply choose $n$, $s$ and $h$ to
balance three errors. In the multi-level algorithms, we choose
$n_\ell,s_\ell,h_\ell$ to minimize the total cost for a fixed total error
using Lagrange multiplier arguments.

For the lognormal case we have similar first order results, see
\cite{GKNSSS15,KSSSU}. There is no higher order results for the lognormal
case because presently there is no QMC theory in this setting.

\section{Software} \label{sec:soft}

The software package QMC4PDE accompanies the survey \cite{KN16}, see
\ifspringer \else \newline \fi
\url{https://people.cs.kuleuven.be/~dirk.nuyens/qmc4pde/}. Here we very
briefly outline its usage.

\subsubsection*{Construction of the generating vector in Python}

In the analysis for the PDE problems we obtain generic bounds on mixed
derivatives of the form
\[
 |\partial^{\bsnu} F(\bsy)| \,\lesssim\,
 \big( (|\bsnu|+a_1)! \big)^{d_1} \prod_{j=1}^s (a_2 \calB\!_j)^{\nu_j}
 \exp (a_3 \calB\!_j |y_j|),
\]
for some constants $a_1$, $a_2$, $a_3$, $d_1$ and some sequence
$\calB\!_j$, where
\[
 F(\bsy) \,=\,
 \begin{cases}
 G(u^s_h) & \mbox{for single-level algorithms}, \\
 G(u^s_{h_\ell} - u^s_{h_{\ell-1}}) & \mbox{for multi-level algorithms},
 \end{cases}
\]
and in particular
\[
 \begin{cases}
 a_3 = 0 & \mbox{for the uniform case}, \\
 a_3 > 0 & \mbox{for the lognormal case}.
 \end{cases}
\]
The Python construction script takes the number of points (as a power
of~2), the dimension, and all these parameters as input from the user,
works out the appropriate weights $\gamma_\setu$, and then constructs a
good generating vector for the QMC rule. This is either a lattice sequence
(constructed following a minimax strategy as described in \cite{CKN06}) or
an interlaced polynomial lattice rule. In the latter case the script also
assembles the interlaced generating matrices, because this is the most
convenient way to generate the points.

\begin{itemize}
\item To construct a generating vector for a lattice sequence (output
    written to file \texttt{z.txt})
\end{itemize}

\begingroup
\ifspringer
\makeatletter
\@totalleftmargin=-0.4cm
\fi
{\small\begin{verbatim}
## uniform case, 100-dim, 2^10 points, with specified bounds b:
./lat-cbc.py --s=100 --m=10 --d2=3 --b="0.1 * j**-3 / log(j+1)"
\end{verbatim}}

{\small\begin{verbatim}
## lognormal case, 100-dim, 2^10 points, with algebraic decay:
./lat-cbc.py --s=100 --m=10 --a2="1/log(2)" --a3=1 --d2=3 --c=0.1
\end{verbatim}}
\endgroup

\begin{itemize}
\item To construct generating matrices for an interlaced polynomial lattice rule
(output written to file \texttt{Bs53.col})
\end{itemize}

\begingroup
\ifspringer
\makeatletter
\@totalleftmargin=-0.4cm
\fi
{\small\begin{verbatim}
## 100-dim, 2^10 points, interlacing 3, with bounds from file:
./polylat-cbc.py --s=100 --m=10 --alpha=3 --a1=5 --b_file=in.txt
\end{verbatim}}
\endgroup

\subsubsection*{Point generators in Matlab/Octave (also available in C++ and Python)}

Here are some Matlab/Octave usage examples for generating the actual QMC
point sets from the output files of the Python construction script.

\begin{itemize}
\item To generate a lattice sequence (specified by the file \texttt{z.txt})
\end{itemize}
{\small\begin{verbatim}
load z.txt                     % load generating vector
latticeseq_b2('init0', z)      % initialize the generator
Pa = latticeseq_b2(20, 512);   % first 512 20-dim points
Pb = latticeseq_b2(20, 512);   % next 512 20-dim points
\end{verbatim}}
\begin{itemize}
\item To generate an interlaced polynomial lattice rule (specified by the file \texttt{Bs53.col})
\end{itemize}
{\small\begin{verbatim}
load Bs53.col                    % load generating matrices
digitalseq_b2g('init0', Bs53)    % initialize the generator
Pa = digitalseq_b2g(100, 512);   % first 512 100-dim points
Pb = digitalseq_b2g(100, 512);   % next 512 100-dim points
\end{verbatim}}

\noindent The same function \texttt{digitalseq\_b2g} can also be used to
generate interlaced Sobol$'$ points by specifying the corresponding
interlaced generating matrices. The parameters for generating Sobol$'$
points are taken from \cite{JK08}.

\begin{itemize}
\item To generate an interlaced Sobol$'$ sequence (interlaced matrices
    specified by the file \texttt{sobol\_alpha3\_Bs53.col})
\end{itemize}
{\small\begin{verbatim}
load sobol_alpha3_Bs53.col           % load generating matrices
digitalseq_b2g('init0', sobol_alpha3_Bs53)   % initialize
Pa = digitalseq_b2g(50, 512);}               % first 512 50-dim
Pb = digitalseq_b2g(50, 512);}               % next 512 50-dim
\end{verbatim}}

\noindent The last example produces interlaced Sobol$'$ points with
interlacing factor $\alpha=3$. They can provide third order convergence if
the integrand has sufficient smoothness.

\section{Concluding Remarks} \label{sec:conc}

QMC (deterministic or randomized) convergence rate and implied constant
can be independent of the dimension. This is achieved by working in a
weighted function space setting. To apply QMC theory, we need an estimate
of the norm of the integrand, and in turn this can help us to choose
appropriate weights for the function space. The chosen weights then enter
the fast CBC construction of the generating vector for the QMC points. The
pairing between the function space setting and the QMC method is very
important, in the sense that we want to achieve the best possible
convergence rate under the weakest assumption on the problem. In practice,
it may be that an off-the-shelf QMC rule works just as well, barring no
theory.

In this article we considered multi-level algorithms. There are other cost
saving strategies for the lognormal case and for other general situations,
see e.g., \cite{DKLS15,HNT16} as well as \cite{GilKS,KN-plenary} in this
volume. Moreover, there have been many others developments on the
application of QMC to PDEs with random coefficients, for some examples see
the last part of the introduction.

%%%%%%%%%%%%%%%%%%%%%%%%%%%%%%%%%%%%%%%%%%%%%%%%%%%%%%%%%%%%%%%%%%%%%%%%%%%%%%%%%%%%%%%%%%%
%%% The acknowledgements
\begin{acknowledgement}
The authors acknowledge the financial supports from the Australian
Research Council (FT130100655 and DP150101770) and the KU Leuven research
fund (OT:3E130287 and C3:3E150478).
\end{acknowledgement}

%%%%%%%%%%%%%%%%%%%%%%%%%%%%%%%%%%%%%%%%%%%%%%%%%%%%%%%%%%%%%%%%%%%%%%%%%%%%%%%%%%%%%%%%%%%
%%% The bibliography
%
% BibTeX users please use
\bibliographystyle{spmpsci}
\bibliography{mybibfile}
% and then copy paste the contents of the .bbl file here for the final version.
%
% E.g.:
%\begin{thebibliography}{99.}%

%\bibitem{ACN2013}
%N.~Achtsis, R.~Cools and D.~Nuyens.
%\newblock Conditional sampling for barrier option pricing under the Heston model.
%\newblock In J.~Dick, F.~Y.\ Kuo, G.~W.\ Peters and I.~H.\ Sloan, editors, {\em {M}onte {C}arlo
%  and Quasi-{M}onte {C}arlo Methods 2012}, pages 253--269. Springer-Verlag, 2013.

%\bibitem{CKN2006}
%R.~Cools, F.~Y. Kuo, and D.~Nuyens.
%\newblock Constructing embedded lattice rules for multivariate integration.
%\newblock {\em SIAM Journal on Scientific Computing}, 28(6):2162--2188, 2006.

%\bibitem{DP2010}
%J.~Dick and F.~Pillichshammer.
%\newblock {\em Digital Nets and Sequences: Discrepancy Theory and Quasi-Monte
%  Carlo Integration}.
%\newblock Cambridge University Press, 2010.

%\bibitem{IT2006}
%J.~{Imai} and K.~S.\ {Tan}.
%\newblock A general dimension reduction technique for derivative pricing.
%\newblock {\em Journal of Computational Finance}, 10(2):129--155, 2006.

%\bibitem{LEC2009}
% P.~L'{\'E}cuyer.
%\newblock Quasi-Monte Carlo methods with applications in finance.
%\newblock {\em Finance and Stochastics}, 13(3):307--349, 2009.

%\end{thebibliography}

\end{document}